\documentclass{amsart}

\usepackage{amsmath,amsthm,amsfonts,amscd,amssymb}
\usepackage{url}

\newtheorem{thm}{Theorem}[section]

\newtheorem{lem}[thm]{Lemma}

\theoremstyle{definition}

\theoremstyle{remark}
\newtheorem{rem}{Remark}[section]

\begin{document}

\title{Bounds on generalized Frobenius numbers}
\author{Lenny Fukshansky \and Achill Sch\"{u}rmann}

\address{Department of Mathematics, 850 Columbia Avenue, Claremont McKenna College, Claremont, CA 91711, USA}
\email{lenny@cmc.edu}
\address{Institute of Mathematics, University of Rostock, 18051 Rostock, Germany}
\email{achill.schuermann@uni-rostock.de}


\begin{abstract}
Let $N \geq2$ and let $1 < a_1 < \dots < a_N$ be relatively prime integers. The Frobenius number of this $N$-tuple is defined to be the largest positive integer that has no representation as $\sum_{i=1}^N a_i x_i$ where $x_1,...,x_N$ are non-negative integers. More generally, the $s$-Frobenius number is defined to be the largest positive integer that has precisely $s$ distinct representations like this. We use techniques from the Geometry of Numbers to give upper and lower bounds on the $s$-Frobenius number for any nonnegative integer $s$.
\end{abstract}

\maketitle

\def\A{{\mathcal A}}
\def\AA{{\mathfrak A}}
\def\B{{\mathcal B}}
\def\C{{\mathcal C}}
\def\D{{\mathcal D}}
\def\EE{{\mathfrak E}}
\def\F{{\mathcal F}}
\def\x{{\mathcal H}}
\def\I{{\mathcal I}}
\def\II{{\mathfrak I}}
\def\J{{\mathcal J}}
\def\K{{\mathcal K}}
\def\kk{{\mathfrak K}}
\def\L{{\mathcal L}}
\def\LL{{\mathfrak L}}
\def\M{{\mathcal M}}
\def\mm{{\mathfrak m}}
\def\MM{{\mathfrak M}}
\def\N{{\mathcal N}}
\def\OO{{\mathfrak O}}
\def\PP{{\mathfrak P}}
\def\R{{\mathcal R}}
\def\PNR{{\mathcal P_N(\real)}}
\def\PMNR{{\mathcal P^M_N(\real)}}
\def\PdNR{{\mathcal P^d_N(\real)}}
\def\s{{\mathcal S}}
\def\V{{\mathcal V}}
\def\X{{\mathcal X}}
\def\Y{{\mathcal Y}}
\def\Z{{\mathcal Z}}
\def\H{{\mathcal H}}
\def\BB{{\mathbb B}}
\def\cee{{\mathbb C}}
\def\pee{{\mathbb P}}
\def\que{{\mathbb Q}}
\def\QQ{{\mathbb Q}}
\def\real{{\mathbb R}}
\def\RR{{\mathbb R}}
\def\zed{{\mathbb Z}}
\def\ZZ{{\mathbb Z}}
\def\aaa{{\mathbb A}}
\def\ff{{\mathbb F}}
\def\VV{{\mathbb V}}
\def\kk{{\mathfrak K}}
\def\qbar{{\overline{\mathbb Q}}}
\def\kbar{{\overline{K}}}
\def\ybar{{\overline{Y}}}
\def\kkbar{{\overline{\mathfrak K}}}
\def\ubar{{\overline{U}}}
\def\eps{{\varepsilon}}
\def\ahat{{\hat \alpha}}
\def\bhat{{\hat \beta}}
\def\gt{{\tilde \gamma}}
\def\h{{\tfrac12}}
\def\be{{\boldsymbol e}}
\def\bei{{\boldsymbol e_i}}
\def\bc{{\boldsymbol c}}
\def\bm{{\boldsymbol m}}
\def\bk{{\boldsymbol k}}
\def\bi{{\boldsymbol i}}
\def\bl{{\boldsymbol l}}
\def\bq{{\boldsymbol q}}
\def\bu{{\boldsymbol u}}
\def\bt{{\boldsymbol t}}
\def\bs{{\boldsymbol s}}
\def\bv{{\boldsymbol v}}
\def\bw{{\boldsymbol w}}
\def\bx{{\boldsymbol x}}
\def\bX{{\boldsymbol X}}
\def\bz{{\boldsymbol z}}
\def\bwy{{\boldsymbol y}}
\def\bY{{\boldsymbol Y}}
\def\bL{{\boldsymbol L}}
\def\ba{{\boldsymbol a}}
\def\baa{{\boldsymbol \alpha}}
\def\bb{{\boldsymbol b}}
\def\bet{{\boldsymbol\eta}}
\def\bxi{{\boldsymbol\xi}}
\def\bo{{\boldsymbol 0}}
\def\bone{{\boldsymbol 1}}
\def\bol{{\boldsymbol 1}_L}
\def\ep{\varepsilon}
\def\p{\boldsymbol\varphi}
\def\q{\boldsymbol\psi}
\def\rank{\operatorname{rank}}
\def\aut{\operatorname{Aut}}
\def\lcm{\operatorname{lcm}}
\def\sgn{\operatorname{sgn}}
\def\spn{\operatorname{span}}
\def\md{\operatorname{mod}}
\def\Norm{\operatorname{Norm}}
\def\dim{\operatorname{dim}}
\def\det{\operatorname{det}}
\def\Vol{\operatorname{Vol}}
\def\rk{\operatorname{rk}}
\def\ord{\operatorname{ord}}
\def\ker{\operatorname{ker}}
\def\div{\operatorname{div}}
\def\Gal{\operatorname{Gal}}
\def\GL{\operatorname{GL}}
\def\SNR{\operatorname{SNR}}
\def\WR{\operatorname{WR}}
\def\scg{\operatorname{\left< \Gamma \right>}}
\def\swrh{\operatorname{Sim_{WR}(\Lambda_h)}}
\def\ch{\operatorname{C_h}}
\def\cht{\operatorname{C_h(\theta)}}
\def\scgt{\operatorname{\left< \Gamma_{\theta} \right>}}
\def\scgmn{\operatorname{\left< \Gamma_{m,n} \right>}}
\def\gat{\operatorname{\Omega_{\theta}}}
\def\GG{\operatorname{G_{\Lambda_{\ba}}}}
\def\Gg{\operatorname{G}}

\section{Introduction}
\label{intro}

Let $N \geq 2$ be an integer and let $a_1,...,a_N$ be positive relatively prime integers. We say that a positive integer $t$ is {\it representable} by the $N$-tuple $\ba := (a_1,...,a_N)$ if
\begin{equation}
\label{rep1}
t = a_1x_1 + \dots + a_Nx_N
\end{equation}
for some nonnegative integers $x_1,\dots,x_N$, and we call each such solution $\bx:=(x_1,\dots,x_N)$ of \eqref{rep1} a {\it representation for $t$ in terms of $\ba$}.  The {\it Frobenius number} $g=g(a_1,...,a_N)$ of this $N$-tuple is defined to be the largest positive integer that has no representations. The condition $\gcd(a_1,...,a_N)=1$ implies that such $g$ exists. More generally, as defined by Beck and Robins in \cite{beck_robins}, let $s$ be a nonnegative integer, and define the {\it $s$-Frobenius number} $g_s=g_s(a_1,...,a_N)$ of $\ba$ to be the largest positive integer that has precisely $s$ distinct representations in terms of $\ba$. Then in particular $g=g_0$.

The Frobenius number has been studied extensively by a variety of authors, starting as early as late 19-th century; see \cite{frobenius_book} for a detailed account and bibliography. More recently, some authors also started studying the more general $s$-Frobenius numbers; for instance, in \cite{shallit} and \cite{beck_kifer} the authors investigated families of $N$-tuples $\ba$ on which the difference $g_s-g_0$ grows unboundedly. This motivates a natural question: how big and how small can $g_s$ be in general?

The main goal of this note is to extend the geometric method of \cite{me_sinai} to obtain general upper and lower bounds on $g_s$.
\smallskip

\begin{rem} \label{other_s_frob} We should warn the reader that the term $s$-Frobenius number is also used by some authors to denote not the largest positive integer that has {\it precisely} $s$ distinct representations in terms of $\ba$, as we do here, but the largest positive integer that has {\it at most} $s$ distinct representations in terms of $\ba$.
\end{rem}

\begin{rem} \label{other_gen} It should also be mentioned that other generalizations of the Frobenius number of different nature have also been considered by a variety of authors. In particular, see Chapter 6 of \cite{frobenius_book}, as well as more recent works \cite{aliev_henk}, \cite{ponomarenko}, and \cite{yoshida},  among others, for further information and references.
\end{rem}
\bigskip

\section{Results}
\label{results}

We start by setting up some notation, following \cite{me_sinai}. Let
$$L_{\ba}(\bX) = \sum_{i=1}^N a_i X_i,$$
be the linear form in $N$ variables with coefficients $a_1,\dots,a_N$, and define the lattice
$$\Lambda_{\ba} = \left\{ \bx \in \zed^N : L_{\ba}(\bx) = 0 \right\}.$$
Let $V_{\ba} = \spn_{\real} \Lambda_{\ba}$, then $V_{\ba}$ is an $(N-1)$-dimensional subspace of $\real^N$ and $\Lambda_{\ba} = V_{\ba} \cap \zed^N$ is a lattice of full rank in $V_{\ba}$. The {\it covering radius} of $\Lambda_{\ba}$ is defined to be
\begin{equation}
\label{cover_radius}
R_{\ba} := \inf \left\{ R \in \real_{>0} : \Lambda_{\ba} + \BB_{V_{\ba}}(R) = V_{\ba} \right\},
\end{equation}
where $\BB_{V_{\ba}}(R)$ is the closed $(N-1)$-dimensional ball of radius $R$ centered at the origin in $V_{\ba}$. For each $1 \leq m \leq N-1$ define the $m$-th {\it successive minimum} of~$\Lambda_{
\ba}$~to~be
\begin{equation}
\label{s_min}
\lambda_m := \min \{ \lambda \in \real : \dim \left( \spn_{\real} \left( \BB_{V_{\ba}}(\lambda) \cap \Lambda_{\ba} \right) \right) \geq m \},
\end{equation}
so $0 < \lambda_1 \leq \dots \leq \lambda_{N-1}$. We also write $\kappa_m$ for the volume of an $m$-dimensional unit ball ($\kappa_0=1$), and $\tau_m$ for the {\it kissing number} in dimension $m$, i.e., the maximal number of unit balls in $\real^m$ that can touch another unit ball. Finally, let us write $\baa_i := (a_1,\dots,a_{i-1},a_{i+1},\dots,a_N)$. We can now state our main results, starting with the upper bounds on $g_s(\ba)$.

\begin{thm} \label{main_upp} With the notation above,
\begin{equation}
\label{mn_bound_u1}
g_s(\ba) \leq \max \left\{ \frac{R_{\ba} (N-1) \sum_{i=1}^N \|\baa_i\| a_i}{\|\ba\|} + 1, \left( s (N-1)! \prod_{i=1}^N a_i \right)^{\frac{1}{N-2}} \right\},
\end{equation}
where $\|\ \|$ stands for the usual Euclidean norm on vectors. If in addition $s \leq \tau_{N-1}+1$, then
\begin{equation}
\label{mn_bound_u2}
g_s(\ba) \leq \frac{3 R_{\ba} \sum_{i=1}^N \|\baa_i\| a_i}{\|\ba\|}.
\end{equation}
\end{thm}
\smallskip

\begin{rem} \label{thm_upp_rem} Note that the quantity $ \frac{R_{\ba} (N-1) \sum_{i=1}^N \|\baa_i\| a_i}{\|\ba\|} + 1$ in the upper bound \eqref{mn_bound_u1} above is precisely the upper bound for the Frobenius number $g_0$ obtained in Theorem 1.1 of \cite{me_sinai}.
\end{rem}

Next we turn to lower bounds. Define the dimensional constant
\begin{equation}
\label{C_N}
C_N = \frac{2^{N^2-\frac{7N}{2}+2} (N-1)^{\frac{N}{2}} ((N-1)!)^{N-1}}{\pi^{\frac{N-2}{2}} \kappa_{N-1}^{N-2}}.
\end{equation}
\begin{thm} \label{main_low} With the notation above,
\begin{equation}
\label{mn_bound_l1}
g_s(\ba) \geq \left( (s+1-N) \prod_{i=1}^N a_i \right)^{\frac{1}{N-1}}.
\end{equation}
Now let $\rho > 1$ be a real number, and suppose that
\begin{equation}
\label{mn_bound_l2_s}
s \geq \frac{\left( \prod_{i=1}^N a_i \right)^{N-2}}{(N-1)!} \left( \frac{C_N \lambda_{N-1}^{N-1}}{\rho-1} \right)^{N-1},
\end{equation}
then
\begin{equation}
\label{mn_bound_l2}
g_s(\ba) \geq \left( \frac{s (N-1)!}{\rho} \prod_{i=1}^N a_i \right)^{\frac{1}{N-1}}.
\end{equation}
\end{thm}
\smallskip

\begin{rem} \label{thm_low_rem} Compare the lower bounds of \eqref{mn_bound_l1} and \eqref{mn_bound_l2} above to the lower bound on the Frobenius number obtained by R\o dseth \cite{rodseth} (see also Theorem 1.1 of \cite{iskander}):
\begin{equation}
\label{frob_lower_bnd}
g_0 \geq \left( (N-1)! {\prod_{i=1}^N a_i} \right)^{\frac{1}{N-1}}.
\end{equation}
In fact, Aliev and Gruber in \cite{iskander} produced a sharp lower bound for $g_0$ in terms of the absolute inhomogeneous minimum of the standard simplex, from which a stronger version of \eqref{frob_lower_bnd} (with a strict inequality) follows. It should also be remarked that the quantities $R_{\ba}$ and $\lambda_{N-1}$, present in our inequalities, can be explicitly bounded using standard techniques from the geometry of numbers. Notice that we can assume without loss of generality that no $a_i$ can be expressed as a nonnegative integer linear combination of the rest of the $a_j$'s: otherwise, $g_s(\ba) = g_s(\baa_i)$. Then equations (28) and (30) of \cite{me_sinai} imply that
\begin{equation}
\label{r_bound}
R_{\ba} \leq \frac{N-1}{2} \lambda_{N-1} \leq \frac{(N-1) \lambda_{N-1}}{\lambda_1} \left( \frac{\|\ba\|}{\kappa_{N-1}} \right)^{\frac{1}{N-1}} \leq  \frac{(N-1) \|\ba\|}{\kappa_{N-1}},
\end{equation}
while equations (25) and (26) of \cite{me_sinai} combined with Minkowski's successive minima theorem (see, for instance, \cite{cass:geom}, p. 203) imply that
\begin{equation}
\label{s_bound}
2 \left( \frac{\|\ba\|}{\kappa_{N-1}(N-1)!} \right)^{\frac{1}{N-1}} \leq \lambda_{N-1} \leq \frac{2 \|\ba\|}{\kappa_{N-1}}.
\end{equation}
In fact, in the situation when the lattice $\Lambda_{\ba}$ is {\it well-rounded} (abbreviated WR), meaning that $\lambda_1 = \dots = \lambda_{N-1}$, inequalities \eqref{r_bound} and \eqref{s_bound} can clearly be improved:
\begin{equation}
\label{wr_bound}
R_{\ba} \leq (N-1) \left( \frac{\|\ba\|}{\kappa_{N-1}} \right)^{\frac{1}{N-1}},\ 2 \left( \frac{\|\ba\|}{\kappa_{N-1}(N-1)!} \right)^{\frac{1}{N-1}} \leq \lambda_{N-1} \leq \left( \frac{2 \|\ba\|}{\kappa_{N-1}} \right)^{\frac{1}{N-1}},
\end{equation}
when $\Lambda_{\ba}$ is WR. The behavior of the Frobenius number $g_0(\ba)$ in this situation was separately studied in \cite{me_sinai}, where WR lattices were called ESM lattices, which stands for {\it equal successive minima}. Finally, the kissing number $\tau_{N-1}$ can be bounded as follows (see pp. 23-24~of~\cite{conway}):
\begin{equation}
\label{tau_bound}
2^{0.2075\dots(N-1)(1+o(1))} \leq \tau_{N-1} \leq 2^{0.401(N-1)(1+o(1))}.
\end{equation}
\end{rem}
\smallskip

We prove Theorems \ref{main_upp} and \ref{main_low} in Section~\ref{proof}. In Section~\ref{count_points} we develop a lattice point counting mechanism, which is used to derive the lower bound of~\eqref{mn_bound_l2}. We are now ready to proceed.
\bigskip

\section{Counting lattice points in polytopes}
\label{count_points}

In this section we present an estimate on the number of lattice points in polytopes, which, while also of independent interest, will be used in Section~\ref{proof} below to prove our main result. To start with, let $P \subset \real^N$ be a polytope of dimension $n \leq N$, i.e., $\dim \VV(P) = n$ where $\VV(P) := \spn_{\real} P$, and let $L \subset \VV(P)$ be a lattice of rank $n$. Define the counting function
$$\Gg(L,P) := \left| L \cap P \right|.$$
Erhart theory studies the properties of $\Gg(L,tP)$ for $t \in \zed_{>0}$, which is a polynomial in $t$ if $P$ is a lattice polytope and a quasipolynomial in $t$ if $P$ is a rational polytope; very little is known in the irrational case (see for instance \cite{beck_robins_book} for a detailed exposition of Erhart theory). In fact, even in the case of a lattice or rational polytope the coefficients of the (quasi-) polynomial $\Gg(L,tP)$ are largely unknown, and hence for many actual applications estimates are needed. Here we record a convenient upper bound on $\Gg(L,P)$. The basic principle going back to Lipschitz (see p. 128 of \cite{lang}) used for such estimates states that when the $n$-dimensional volume $\Vol_n(P)$ is large comparing to $\det(L)$, then $\Gg(L,P)$ can be approximated by $\frac{\Vol_n(P)}{\det(L)}$, and so the problem comes down to estimating the error term of such approximation. An upper bound on this error term -- not only for polytopes, but for a rather general class of compact domains -- has been produced by Davenport \cite{davenport} and then further refined by Thunder \cite{thunder}. Here we present a variation of Thunder's bound in case of polytopes. 

Generalizing the notation of Section~\ref{intro} to arbitrary lattices, let $\BB_{\VV(P)}(R)$ be a ball of radius $R$ centered at the origin in $\VV(P)$, and for each $1 \leq m \leq n$ define the $m$-th successive minimum of $L$ as in \eqref{s_min} above:
$$\lambda_m = \min \{ \lambda \in \real : \dim \left( \spn_{\real} \left( \BB_{\VV(P)}(\lambda) \cap L \right) \right) \geq m \}.$$
Also for each $1 \leq m \leq n$, let
\begin{equation}
\label{V_m}
V_m(P) := \max \{ \Vol_m(F) : F \text{ is an $m$-dimensional face  of } P \}.
\end{equation}
With this notation at hand, the following estimate is an immediate implication of Theorem~4 of \cite{thunder}.

\begin{lem} \label{count_lemma} With notation as above,
$$\Gg(L,P) \leq \frac{\Vol_n(P)}{\det(L)} + \sum_{m=0}^{n-1} \frac{2^{(n+1)m} \left( m n! \right)^m}{\kappa_m \kappa_n^m} \binom{n}{m} \frac{V_m(P)}{\lambda_1 \cdots \lambda_m},$$
where the product $\lambda_1 \dots \lambda_m$ is interpreted as 1 when $m=0$.
\end{lem}
\smallskip

\begin{rem} \label{count_polytope} Notice that Lemma~\ref{count_lemma}, and more generally the counting estimates discussed in section~5 of \cite{thunder}, provide a mechanism for producing explicit polynomial bounds on the number of points of an arbitrary lattice in a variety of homogeneously expanding compact domains, which is especially easy to use in case of polytopes (as we do in Section~\ref{proof} for certain simplices). This observation gives a partial solution to Problem~3.2 of \cite{snowbird}, previously formulated by the first author.
\end{rem}
\smallskip

In the next section, we apply Lemma~\ref{count_lemma} to derive the lower bound~of~\eqref{mn_bound_l2}.
\bigskip

\section{Bounds on $g_s(\ba)$}
\label{proof}

In this section we prove Theorems \ref{main_upp} and \ref{main_low}, deriving the inequalities \eqref{mn_bound_u1}, \eqref{mn_bound_u2}, \eqref{mn_bound_l1}, and \eqref{mn_bound_l2}. For a positive integer $t$, consider the hyperplane $V_{\ba}(t)$ in $\real^N$ defined by the equation \eqref{rep1}, which is a translate of $V_{\ba}$, and write $\Lambda_{\ba}(t) = V_{\ba}(t) \cap \zed^N$. Fix a point $\bu_t \in \Lambda_{\ba}(t)$, and define a translation map $f_t : V_{\ba} \rightarrow V_{\ba}(t)$ given by $f_t(\bx) = \bx + \bu_t$ for each $\bx \in V_{\ba}$. Then $f_t$ is bijective and preserves distance; moreover, it maps $\Lambda_{\ba}$ bijectively onto $\Lambda_{\ba}(t)$. The intersection of $V_{\ba}(t)$ with the positive orthant $\real^N_{\geq 0}$ is an $(N-1)$-dimensonal simplex, call it $S(t)$. Then define
\begin{equation}
\label{GG}
\GG(t) := | \Lambda_{\ba}(t) \cap S(t) | = | \Lambda_{\ba} \cap f_t^{-1}(S(t)) |,
\end{equation}
and notice that each point in $\Lambda_{\ba}(t) \cap S(t)$ corresponds to a solution of \eqref{rep1} in non-negative integers. Hence for every $t > g_s(\ba)$ we have $\GG(t) > s$. Moreover, $g_s(\ba)$ is precisely the smallest among all positive integers $m$ such that for each integer $t > m$, $\GG(t) > s$. Therefore, in order to obtain bounds on $g_s(\ba)$, we want to produce estimates on $\GG(t)$, which is what we do next.
\smallskip

Combining \eqref{GG} with bounds by Blichfeldt \cite{blichfeldt} (see also equation (3.2) of \cite{lattice_points}) and by Gritzmann \cite{gritzmann} (see also equation (3.3) of \cite{lattice_points}), we have:
\begin{equation}
\label{lattice_bounds_1}
\frac{\Vol_{N-1}(S(t)) - R_{\ba} A_{N-1}(S(t))}{\det \Lambda_{\ba}} \leq \GG(t) \leq \frac{\Vol_{N-1}(S(t))}{\det \Lambda_{\ba}} (N-1)! + (N-1),
\end{equation}
where $\Vol_{N-1}(S(t))$ is the volume and $A_{N-1}(S(t))$ is the surface area of $S(t)$, and $R_{\ba}$ is the covering radius of $\Lambda_{\ba}$ as defined in \eqref{cover_radius} above. Equations (17) and (18) of \cite{me_sinai} state that
\begin{equation}
\label{vol_area}
\Vol_{N-1}(S(t)) = \frac{t^{N-1} \|\ba\|}{(N-1)! \prod_{i=1}^N a_i},\ \ A_{N-1}(S(t)) = \frac{t^{N-2} \sum_{i=1}^N \|\baa_i\| a_i}{(N-2)! \prod_{i=1}^N a_i}.
\end{equation}
In addition, by equation (25) of \cite{me_sinai}, $\det \Lambda_{\ba} = \|\ba\|$. Combining these observations with \eqref{lattice_bounds_1}, we obtain
\begin{equation}
\label{lattice_bounds_2.2}
\GG(t) \geq \frac{t^{N-2}}{(N-2)! \prod_{i=1}^N a_i} \left( \frac{t}{N-1} - \frac{R_{\ba} \sum_{i=1}^N \|\baa_i\| a_i}{\|\ba\|} \right),
\end{equation}
and
\begin{equation}
\label{lattice_bounds_2.1}
\GG(t) \leq \frac{t^{N-1}}{\prod_{i=1}^N a_i} + (N-1).
\end{equation}
\smallskip

Notice however that Blichfeldt's upper bound of \eqref{lattice_bounds_1} is weaker than the bound of Lemma \ref{count_lemma} for large $t$, hence our next goal is to produce an explicit upper bound on $\GG(t)$ from Lemma \ref{count_lemma}. Since each $m$-dimensional face of $S(t)$ is an $m$-dimensional simplex for each $0 \leq m \leq N-1$, equation (17) of \cite{me_sinai} implies that
\begin{equation}
\label{vol_sf}
V_m(S(t)) \leq \frac{t^m \|\ba\|}{m!}.
\end{equation}
On the other hand, Minkowski's successive minima theorem implies that for each $1 \leq m \leq N-2$,
\begin{equation}
\label{suc_min_b}
\lambda_1 \dots \lambda_m \geq \frac{2^{N-1} \det \Lambda_{\ba}}{(N-1)! \lambda_{m+1} \dots \lambda_{N-1}} \geq \frac{2^{N-1} \|\ba\|}{(N-1)! \lambda_{N-1}^{N-1-m}}.
\end{equation}
Also notice that for all $1 \leq m \leq N-1$,
\begin{equation}
\label{fact_ineq}
\frac{m^m}{\kappa_m m!} = \frac{m^m \Gamma \left( 1+\frac{m}{2} \right)}{\pi^{m/2}m!} = 
\left\{ \begin{array}{ll}
\frac{(2k)^{2k} k!}{\pi^{k} (2k)!} & \mbox{if $m=2k$} \\
\frac{(2k+1)^{2k+1}}{\pi^{k} 2^{2k+1} k!} & \mbox{if $m=2k+1$}
\end{array}
\right. \leq \left( \frac{2m}{\pi} \right)^{m/2},
\end{equation}
where $\Gamma$ stands for the $\Gamma$-function. Finally, $\binom{N-1}{m} \leq (N-1) \binom{N-2}{m}$. Define
\begin{equation}
\label{C1N}
C'_N = \frac{(N-1) (N-1)!}{2^{N-1}}.
\end{equation}
Combining \eqref{vol_sf}, \eqref{suc_min_b}, and \eqref{fact_ineq} with Lemma \ref{count_lemma}, we obtain:
\begin{eqnarray}
\label{lattice_bounds_thunder}
\GG(t) & \leq &  \frac{t^{N-1}}{(N-1)! \prod_{i=1}^N a_i} + C'_N \lambda_{N-1} \times  \nonumber \\
& & \times\ \sum_{m=0}^{N-2} \binom{N-2}{m} \left( \frac{2^N (N-2)^{1/2} (N-1)!\ t}{\kappa_{N-1} \sqrt{2\pi}} \right)^m   \lambda_{N-1}^{N-2-m} \nonumber \\
& \leq & \frac{t^{N-1}}{(N-1)! \prod_{i=1}^N a_i} + C'_N \lambda_{N-1} \left(\frac{2^N (N-2)^{1/2} (N-1)!\ t}{\kappa_{N-1} \sqrt{2\pi}} +  \lambda_{N-1} \right)^{N-2} \nonumber \\
& \leq & \frac{t^{N-1}}{(N-1)! \prod_{i=1}^N a_i} + \frac{2^{N^2-\frac{7N}{2}+2} (N-1)^{\frac{N}{2}} ((N-1)!\ \lambda_{N-1})^{N-1} t^{N-2}}{\pi^{\frac{N-2}{2}} \kappa_{N-1}^{N-2}}.
\end{eqnarray}
Then for any $\rho > 1$,
\begin{equation}
\label{lb_thunder}
\GG(t) \leq \frac{\rho t^{N-1}}{(N-1)! \prod_{i=1}^N a_i}, \text{ when } t \geq \frac{C_N \lambda_{N-1}^{N-1} \prod_{i=1}^N a_i}{\rho-1},
\end{equation}
where $C_N$ is as in \eqref{C_N}.

\begin{rem} \label{WR_rem} Similarly to the observations in Remark~\ref{thm_low_rem}, the inequality \eqref{suc_min_b} can be improved in case $\Lambda_{\ba}$ is WR. As a result in this case, inequalities \eqref{lattice_bounds_thunder} and \eqref{lb_thunder} can also be made stronger.
\end{rem}
\smallskip

A different technique can be used to produce a lower bound on $\GG(t)$ for small $t$. Notice that an open ball of radius $R_{\ba}$ in $V_{\ba}$ contains at least one point of $\Lambda_{\ba}$, hence one can estimate the number of such balls in $S(t)$ to obtain a lower bound on $\GG(t)$. The kissing number $\tau_{N-1}$ is the maximal number of balls of radius $R_{\ba}$ that can touch another ball of radius $R_{\ba}$ without overlap, hence each ball of radius $3R_{\ba}$ in $V_{\ba}$ contains an arrangement of $\tau_{N-1}+1$ non-overlapping balls of radius $R_{\ba}$. Now a standard isoperimetric identity (see, for instance, equation (1.3) of \cite{betke_henk}) implies that the inradius $r(t)$ of the simplex $S(t)$ satisfies
\begin{equation}
\label{inradius}
r(t) = \frac{(N-1) \Vol_{N-1}(S(t))}{A_{N-1}(S(t))} = \frac{t \|\ba\|}{\sum_{i=1}^N \|\baa_i\| a_i},
\end{equation}
and so if $t \geq  \frac{3R_{\ba} \sum_{i=1}^N \|\baa_i\| a_i}{\|\ba\|}$, then $S(t)$ contains a ball of radius $3R_{\ba}$, and hence at least $\tau_{N-1}+1$ points of $\Lambda_{\ba}$. In other words,
\begin{equation}
\label{t1_inradius}
\GG(t) \geq \tau_{N-1}+1, \text{ when } t \geq  \frac{3R_{\ba} \sum_{i=1}^N \|\baa_i\| a_i}{\|\ba\|}.
\end{equation}
Now, equipped with these inequalities on $\GG(t)$, we can easily derive the bounds of Theorems \ref{main_upp} and \ref{main_low}.
\bigskip

First notice that if we pick $t$ greater than the maximal expression in the upper bound of \eqref{mn_bound_u1}, then \eqref{lattice_bounds_2.2} implies $\GG(t) > s$. In addition, \eqref{t1_inradius} implies that for $s \leq \tau_{N-1}+1$, $g_s(\ba)$ satisfies \eqref{mn_bound_u2}. As for lower bounds on $g_s(\ba)$, if we pick 
$$t \leq \left( (s+1-N) {\prod_{i=1}^N a_i} \right)^{\frac{1}{N-1}},$$
then \eqref{lattice_bounds_2.1} implies $\GG(t) \leq s$, and so produces the lower bound of \eqref{mn_bound_l1}. Finally, \eqref{lb_thunder} implies that when $s$ satisfies \eqref{mn_bound_l2_s}, $g_s(\ba)$ satisfies \eqref{mn_bound_l2}. This completes the proof of Theorems \ref{main_upp} and \ref{main_low}.
\qed
\smallskip

\begin{rem} \label{other_bounds} For comparison purposes with \eqref{lattice_bounds_thunder}, we mention another upper bound on $\GG(t)$, which is given by equation (3.3) of \cite{lattice_points}:
\begin{equation}
\label{lattice_bounds_2.3}
\GG(t) \leq \frac{\Vol_{N-1}(S(t)+C(\Lambda_{\ba}))}{\det \Lambda_{\ba}} \leq \frac{\Vol_{N-1}(S(t)+\BB_{N-1}(R_{\ba}))}{\det \Lambda_{\ba}},
\end{equation}
where
\begin{equation}
\label{voronoi}
C(\Lambda_{\ba}) := \{ \bwy \in V_{\ba} : \| \bwy \| \leq \| \bwy - \bx \|\ \forall\ \bx \in \Lambda_{\ba} \}
\end{equation}
is the Voronoi cell of the lattice $\Lambda_{\ba}$. Now the right hand side of \eqref{lattice_bounds_2.3} can be expanded using mixed volumes (see for instance \cite{mixed_volumes}), i.e.:
\begin{equation}
\label{lattice_bounds_2.4}
\Vol_{N-1}(S(t)+\BB_{N-1}(R_{\ba})) = \sum_{m=0}^{N-1} \kappa_m R_{\ba}^m \V_{N-m-1}(S(t)),
\end{equation}
where $\V_k(S(t))$ denotes the $k$-th mixed volume of $S(t)$. In particular, 
$$\V_{N-1}(S(t)) = \Vol_{N-1}(S(t)),\ \V_{N-2}(S(t)) = \frac{1}{2} A_{N-1}(S(t)),$$
as given by \eqref{vol_area}, and $\V_0(K)=1$. Then combining \eqref{lattice_bounds_2.3}, \eqref{lattice_bounds_2.4}, and \eqref{vol_area}, we obtain an upper bound on $\GG(t)$ in terms of the covering radius $R_{\ba}$, analogous to the lower bound of \eqref{lattice_bounds_2.2}:
\begin{equation}
\label{lattice_bounds_2.5}
\GG(t) \leq \frac{t^{N-1}}{(N-1)! \prod_{i=1}^N a_i} + \frac{t^{N-2} R_{\ba} \sum_{i=1}^N \|\baa_i\| a_i}{(N-2)! \|\ba\| \prod_{i=1}^N a_i} + \sum_{m=2}^{N-1} \frac{\kappa_m  R_{\ba}^m \V_{N-m-1}(S(t))}{\|\ba\|}.
\end{equation}
The bound of \eqref{lattice_bounds_2.5} is similar in spirit to that of \eqref{lattice_bounds_thunder}, although the mixed volumes may generally be hard to compute. An expansion similar to \eqref{lattice_bounds_2.4} has recently been used by M. Henk and J. M. Wills to obtain a strengthening of Blichfeldt's upper bound as in \eqref{lattice_bounds_1}, at least in the case of the integer lattice $\zed^N$ (see Theorem 1.1 and Conjecture 1.1 of \cite{henk_wills}).
\end{rem}
\bigskip

{\bf Acknowledgment.} We would like to thank the anonymous referees for their helpful comments on the subject of this paper.
\bigskip

\section{Appendix: erratum and addendum}
\label{appendix}

Here we correct two inaccuracies in the statement of Theorem 2.2 of the published version of our paper. We also  exhibit additional  bounds on the generalized Frobenius numbers, which complement those developed in the paper.

\subsection{Correction to Theorem 2.2}
\label{erratum}

Let the notation be as above. The proof of formula (7) in Theorem 2.2 depends on Blichfeldt's inequality (20), which is true with the additional assumption that the simplex $f_t^{-1}(S(t))$ contains $N-1$ linearly independent points of $\Lambda_{\ba}$. If $r(t)$, the inradius of $S(t)$, is $\geq \lambda_{N-1}$, the last successive minimum of $\Lambda_{\ba}$, then this condition is satisfied. Using identity (27) for $r(t)$ we can easily deduce that this happens when
\begin{equation}
\label{t1}
t \geq t_* := \frac{\lambda_{N-1} \sum_{i=1}^N \|\baa_i\| a_i}{\|\ba\|}.
\end{equation}
Now assume that \eqref{t1} is not satisfied, i.e. $t < t_*$. In this case, 
$$f_t^{-1}(S(t)) \subset f_{t_*}^{-1}(S(t_*)),$$
and so $\GG(t) \leq \GG(t_*)$, and we can apply Blichfeldt's bound on $\GG(t_*)$. Then we obtain
$$\GG(t) \leq \GG(t_*) \leq \frac{t_*^{N-1}}{\prod_{i=1}^N a_i} + (N-1).$$
In other words, the inequality
\begin{equation}
\label{t2}
\GG(t) \leq \frac{\max \{t,t_*\}^{N-1}}{\prod_{i=1}^N a_i} + (N-1)
\end{equation}
holds for any $t$. Using this inequality instead of (20), we see that if
\begin{equation}
\label{t3}
\max\{ t,t_*\} \leq \left( (s+1-N) \prod_{i=1}^N a_i \right)^{\frac{1}{N-1}},
\end{equation}
then $\GG(t) \leq s$. Now \eqref{t3} holds when
$$s \geq \frac{t_*^{N-1}}{\prod_{i=1}^N a_i} + (N-1).$$
This means that the following addition to the statement of Theorem 2.2 should be made: formula (7) holds under the assumption that
$$s \geq \frac{ \left( \lambda_{N-1} \sum_{i=1}^N \|\baa_i\| a_i\right)^{N-1}}{\|\ba\|^{N-1} \prod_{i=1}^N a_i} + (N-1).$$

\noindent
{\bf Acknowledgment.} We would like to thank Iskander Aliev and Martin Henk for attracting our attention to this inaccuracy.
\smallskip

The second inaccuracy in our Theorem 2.2 comes from the application of J. Thunder's Theorem 4 of \cite{thunder} as recorded in Lemma 3.1 above: in the statement of this lemma, the quantity 
$$V_m(P) := \max \{ \Vol_m(F) : F \text{ is an $m$-dimensional face  of } P \},$$
as defined in equation (15) above should be replaced by
\begin{equation}
\label{v_prime}
V'_m(P) := \sum \Vol_m(F),
\end{equation}
for all $1 \leq m \leq n$, where the sum is over all $m$-dimensional faces of $P$; here $P$ is an $n$-dimensional polytope in $\real^N$, $n \leq N$. Hence Lemma 3.1 should read as follows.

\begin{lem} [Lemma 3.1, corrected] \label{count_lemma_1} With notation as above,
$$\Gg(L,P) \leq \frac{\Vol_n(P)}{\det(L)} + \sum_{m=0}^{n-1} \frac{2^{(n+1)m} \left( m n! \right)^m}{\kappa_m \kappa_n^m} \binom{n}{m} \frac{V'_m(P)}{\lambda_1 \cdots \lambda_m},$$
where $V'_m(P)$ is as in \eqref{v_prime} above and the product of successive minima $\lambda_1 \dots \lambda_m$ is interpreted as 1 when $m=0$.
\end{lem}
\smallskip

\noindent
Lemma 3.1 is applied in case $P$ is the simplex $S(t)$ in the proof of inequality (9) of Theorem 2.2. We can now correct this argument by applying our Lemma~\ref{count_lemma_1} instead. The total number of $m$-faces of $S(t)$ is $\binom{N}{m+1}$ for each $0 \leq m \leq N-1$, and so by formula (21),
\begin{equation}
\label{v_prime_1}
V'_m(S(t)) \leq \binom{N}{m+1} V_m(S(t)) \leq  \binom{N}{m+1} \frac{t^m \|\ba\|}{m!} \leq \frac{N^N t^m \|\ba\|}{N!\ m!}.
\end{equation}
Now we can proceed with the derivation of the inequality (25), applying Lemma~\ref{count_lemma_1} above instead of Lemma 3.1 and inequality \eqref{v_prime_1} above instead of inequality (21), we readily obtain
$$\GG(t) \leq \frac{t^{N-1}}{(N-1)! \prod_{i=1}^N a_i} + \frac{2^{N^2-\frac{7N}{2}+2} N^N (N-1)^{\frac{N}{2}} ((N-1)!\ \lambda_{N-1})^{N-1} t^{N-2}}{N!\ \pi^{\frac{N-2}{2}} \kappa_{N-1}^{N-2}}.$$
Taking all these remarks into account, the correct statement of Theorem 2.2 should be as follows.

\begin{thm} [Theorem 2.2, corrected]  \label{thm_lower} With the notation as above,
\begin{equation}
\label{mn_bound_l1_1}
g_s(\ba) \geq \left( (s+1-N) \prod_{i=1}^N a_i \right)^{\frac{1}{N-1}}
\end{equation}
for all
\begin{equation}
\label{t4}
s \geq \frac{ \left( \lambda_{N-1} \sum_{i=1}^N \|\baa_i\| a_i\right)^{N-1}}{\|\ba\|^{N-1} \prod_{i=1}^N a_i} + (N-1).
\end{equation}
Now let $\rho > 1$ be a real number, and suppose that
\begin{equation}
\label{mn_bound_l2_s_1}
s \geq \frac{\left( \prod_{i=1}^N a_i \right)^{N-2}}{(N-1)!} \left( \frac{A_N \lambda_{N-1}^{N-1}}{\rho-1} \right)^{N-1},
\end{equation}
where 
\begin{equation}
\label{A_N}
A_N = \frac{2^{N^2-\frac{7N}{2}+2} N^N (N-1)^{\frac{N}{2}} ((N-1)!)^{N-1}}{N!\ \pi^{\frac{N-2}{2}} \kappa_{N-1}^{N-2}}.
\end{equation}
Then
\begin{equation}
\label{mn_bound_l2_1}
g_s(\ba) \geq \left( \frac{s (N-1)!}{\rho} \prod_{i=1}^N a_i \right)^{\frac{1}{N-1}}.
\end{equation}
\end{thm}
\smallskip

\bigskip

\subsection{Additional bounds on $s$-Frobenius numbers}
\label{addendum}

Here we exhibit additional bounds on $s$-Frobenius numbers, following the same principle as above. These bounds are of the same order of magnitude as in the theorems above, but may be more convenient in some applications. Let all the notation be as above. 

We first produce a lower bound on $g_s(\ba)$ employing a new lattice point counting estimate due to M. Widmer \cite{widmer}. 

\begin{thm} \label{lower_bnd} With the notation as above,
\begin{equation}
\label{mn_bound_l3}
g_s(\ba) \geq \frac{\left( (N-2)! \right)^{\frac{1}{N-1}}}{4 (N-1)^{\frac{3(N+1)}{2}}} \left( \frac{s R_{\ba} \prod_{i=1}^N a_i}{\|\ba\|} \right)^{\frac{1}{N-1}}
\end{equation}
for all
\begin{equation}
\label{s_rest}
s \geq \frac{\left( 4 \lambda_{N-1} \sum_{i=1}^N \|\baa_i\| a_i\right)^{N-1} (N-1)^{\frac{3N^2-1}{2}}}{(N-1)! \ \|\ba\|^{N-2} R_{\ba} \prod_{i=1}^N a_i}.
\end{equation}
\end{thm}

\proof
To obtain the lower bound \eqref{mn_bound_l3}, we argue in precisely the same way as in Section~4 above, replacing Blichfeldt's upper bound on $\GG(t)$ with the bound of Proposition 2.9 of \cite{widmer}:
\begin{equation}
\label{lattice_bounds_1_1}
\GG(t) \leq 8^{N-1} (N-1)^{\frac{3(N^2-1)}{2}} \frac{\Vol_{N-1}(S(t))}{\lambda_1 \cdots \lambda_{N-1}}.
\end{equation}
This inequality holds under the same assumption as Blichfeldt's bound, namely whenever the simplex $f_t^{-1}(S(t))$ contains $N-1$ linearly independent points of $\Lambda_{\ba}$. This means that $t$ needs to satisfy condition \eqref{t1}, as in Section~\ref{erratum} above, for us to apply \eqref{lattice_bounds_1_1}. Equation (17) of \cite{me_sinai} states that
\begin{equation}
\label{vol}
\Vol_{N-1}(S(t)) = \frac{t^{N-1} \|\ba\|}{(N-1)! \prod_{i=1}^N a_i}.
\end{equation}
In addition, by equation (25) of \cite{me_sinai}, $\det \Lambda_{\ba} = \|\ba\|$. Combining these observations with \eqref{lattice_bounds_1_1} and (26), (28) of \cite{me_sinai}, we obtain
\begin{equation}
\label{lattice_bounds_2.2_1}
\GG(t) \leq \frac{4^{N-1} (N-1)^{\frac{3N^2-1}{2}}}{(N-1)!} \times \frac{\max \{t,t_*\}^{N-1} \|\ba\|}{R_{\ba} \prod_{i=1}^N a_i},
\end{equation}
where $t_*$ is as in \eqref{t1}. Now notice that if 
$$\max \{t,t_*\} \leq \left( \frac{s (N-1)! R_{\ba} \prod_{i=1}^N a_i}{4^{N-1} (N-1)^{\frac{3N^2-1}{2}} \|\ba\|} \right)^{\frac{1}{N-1}},$$
then \eqref{lattice_bounds_2.2_1} implies $\GG(t) \leq s$, and so produces the lower bound of \eqref{mn_bound_l3}. This means, however, that $s$ needs to satisfy
$$s \geq \frac{(4t_*)^{N-1} (N-1)^{\frac{3N^2-1}{2}} \|\ba\|}{(N-1)! R_{\ba} \prod_{i=1}^N a_i},$$
which, combined with \eqref{t1}, produces \eqref{s_rest}.
\endproof
\smallskip

We now produce an upper bound on $g_s(\ba)$ in terms of the ratio of the covering radius $R_{\ba}$ of the lattice $\Lambda_{\ba}$ and the inradius $r(1)$ of the simplex $S(1)$; also notice the exponent $\frac{1}{N-1}$ in this upper bound, which is the same as in all our lower bounds. This is just a variation on the results above.

\begin{thm} \label{upper_bnd} With the notation as above,
\begin{equation}
\label{mn_bound_u2_1}
g_s(\ba) \leq \max \left\{ \frac{(N-1) R_{\ba}}{r(1)} + 1, \left( \left( \frac{(N-1) R_{\ba}}{r(1)} + 1 \right) s (N-1)! \prod_{i=1}^N a_i \right)^{\frac{1}{N-1}} \right\}.
\end{equation}
\end{thm}

\proof
For convenience, define
\begin{equation}
\label{beta}
\beta(\ba) := \frac{R_{\ba} (N-1) \sum_{i=1}^N \|\baa_i\| a_i}{\|\ba\|} + 1.
\end{equation}
Now, equation (19) gives
\begin{equation}
\label{lattice_bounds_2.3_1}
\GG(t) \geq \frac{t^{N-1}}{(N-1)! \prod_{i=1}^N a_i} \left( 1 - \frac{(N-1) R_{\ba} \sum_{i=1}^N \|\baa_i\| a_i}{t \|\ba\|} \right).
\end{equation}
If we assume that $t \geq \beta(\ba)$, then we obtain
\begin{eqnarray}
\label{lattice_bounds_2.4_1}
\GG(t) & \geq & \frac{t^{N-1}}{(N-1)! \prod_{i=1}^N a_i} \left( 1 - \frac{(N-1) R_{\ba} \sum_{i=1}^N \|\baa_i\| a_i}{\beta(\ba) \|\ba\|} \right) \nonumber \\
& = & \frac{t^{N-1}}{(N-1)! \beta(\ba) \prod_{i=1}^N a_i}.
\end{eqnarray}
Now notice that if we pick 
$$t \geq \left( \beta(\ba) s (N-1)! \prod_{i=1}^N a_i \right)^{\frac{1}{N-1}},$$
then \eqref{lattice_bounds_2.4_1} implies $\GG(t) \geq s$. Combining this with the fact that $t$ has to be at least $\beta(\ba)$ produces the following upper bound:
\begin{equation}
\label{mn_bound_u1_1}
g_s(\ba) \leq \max \left\{ \beta(\ba), \left( \beta(\ba) s (N-1)! \prod_{i=1}^N a_i \right)^{\frac{1}{N-1}} \right\}.
\end{equation}
The isoperimetric identity \eqref{inradius} on the inradius $r(t)$ of the simplex $S(t)$ implies that
$$\frac{\sum_{i=1}^N \|\baa_i\| a_i}{\|\ba\|} = \frac{1}{r(1)}.$$
Then the definition of $\beta(\ba)$ implies that
$$\beta(\ba) = \frac{(N-1) R_{\ba}}{r(1)} + 1.$$
This allows us to rewrite the upper bound of \eqref{mn_bound_u1_1} in terms of the ratio of $R_{\ba}$, the covering radius of $\Lambda_{\ba}$, and $r(1)$, the inradius of the simplex $S(1)$, producing \eqref{mn_bound_u2_1}. This completes the proof of Theorem~\ref{upper_bnd}.
\endproof
\bigskip

\begin{rem} Additional bounds of comparable order of magnitude on $s$-Frobenius numbers  have also been produced in \cite{aliev_henk_me} with the use of a rather different method, where they are applied to study the average behavior of the $s$-Frobenius numbers.
\end{rem}
\bigskip

\bibliographystyle{plain}  
\bibliography{frob_rep}    

\begin{thebibliography}{10}

\bibitem{frobenius_book}
J.~L.~Ramirez Alfonsin.
\newblock {\em The Diophantine Frobenius problem}.
\newblock Oxford University Press, 2005.

\bibitem{aliev_henk_me}
I.~Aliev, L.~Fukshansky, and M.~Henk.
\newblock Generalized {F}robenius numbers: bounds and average behavior.
\newblock {\em Acta Arithm., to appear}.

\bibitem{iskander}
I.~Aliev and P.~M. Gruber.
\newblock An optimal lower bound for the {F}robenius problem.
\newblock {\em J. Number Theory}, 123(1):71--79, 2007.

\bibitem{aliev_henk}
I.~Aliev and M.~Henk.
\newblock On feasibility of integer knapsacks.
\newblock {\em SIAM J. Optimization, to appear; arXiv:0911.4186}, 2010.

\bibitem{ponomarenko}
J.~Amos, I.~Pascu, V.~Ponomarenko, E.~Trevino, and Y.~Zhang.
\newblock The multi-dimensional {F}robenius problem.
\newblock {\em Adv. in Appl. Math., to appear}, 2006.

\bibitem{snowbird}
M.~Beck, B.~Chen, L.~Fukshansky, C.~Haase, A.~Knutson, B.~Reznick, S.~Robins,
  and A.~Sch{\"u}rmann.
\newblock Problems from the {C}ottonwood {R}oom.
\newblock In {\em Integer points in polyhedra---geometry, number theory,
  algebra, optimization}, Contemp. Math., 374, pages 179--191. Amer. Math.
  Soc., Providence, RI, 2005.

\bibitem{beck_kifer}
M.~Beck and C.~Kifer.
\newblock An extreme family of generalized {F}robenius numbers.
\newblock {\em preprint; ar{X}iv:1005.2692}, 2010.

\bibitem{beck_robins}
M.~Beck and S.~Robins.
\newblock A formula related to the {F}robenius problem in two dimensions.
\newblock In {\em Number Theory (New York Seminar 2003)}, pages 17--23.
  Springer, New York, 2004.

\bibitem{beck_robins_book}
M.~Beck and S.~Robins.
\newblock {\em Computing the continuous discretely. Integer-point enumeration
  in polyhedra}.
\newblock Springer, 2007.

\bibitem{betke_henk}
U.~Betke and M.~Henk.
\newblock A generalization of {S}teinhagen's theorem.
\newblock {\em Abh. Math. Sem. Univ. Hamburg}, 63:165--176, 1993.

\bibitem{blichfeldt}
H.~F. Blichfeldt.
\newblock Note on the geometry of numbers.
\newblock {\em Bull. Amer. Math. Soc.}, 27:150--153, 1921.

\bibitem{cass:geom}
J.~W.~S. Cassels.
\newblock {\em An Introduction to the Geometry of Numbers}.
\newblock Springer-Verlag, 1959.

\bibitem{conway}
J.~H. Conway and N.~J.~A. Sloane.
\newblock {\em Sphere Packings, Lattices, and Groups}.
\newblock Springer-Verlag, {T}hird edition, 1999.

\bibitem{davenport}
H.~Davenport.
\newblock On a principle of {L}ipschitz.
\newblock {\em J. London Math. Soc.}, 26:179--183, 1951.

\bibitem{rodseth}
\O. J.~R\o dseth.
\newblock An upper bound for the $h$-range of the postage stamp problem.
\newblock {\em Acta Arith.}, 54(4):301--306, 1990.

\bibitem{me_sinai}
L.~Fukshansky and S.~Robins.
\newblock Frobenius problem and the covering radius of a lattice.
\newblock {\em Discrete Comput. Geom.}, 37(3):471--483, 2007.

\bibitem{gritzmann}
P.~Gritzmann.
\newblock {\em Finite {P}ackungen und {\"U}berdeckungen}.
\newblock {H}abiltationsschrift, Universit{\"a}t Siegen, 1984.

\bibitem{lattice_points}
P.~Gritzmann and J.~M. Wills.
\newblock Lattice points.
\newblock In {\em Handbook of Convex Geometry, Vol. A, B}, pages 765--797.
  North-Holland, Amsterdam, 1993.

\bibitem{henk_wills}
M.~Henk and J.~M. Wills.
\newblock A {B}lichfeldt-type inequality for the surface area.
\newblock {\em Monatsh. Math.}, 154:135--144, 2008.

\bibitem{lang}
S.~Lang.
\newblock {\em Algebraic Number Theory}.
\newblock Springer-Verlag, 1986.

\bibitem{mixed_volumes}
J.~R. Sangwine-Yager.
\newblock Mixed volumes.
\newblock In {\em Handbook of Convex Geometry, Vol. A, B}, pages 43--71.
  North-Holland, Amsterdam, 1993.

\bibitem{shallit}
J.~Shallit and J.~Stankewicz.
\newblock Unbounded discrepancy in {F}robenius numbers.
\newblock {\em preprint; ar{X}iv:1003.0021}, 2010.

\bibitem{yoshida}
A.~Takemura and R.~Yoshida.
\newblock A generalization of the integer linear infeasibility problem.
\newblock {\em Discrete Optim.}, 5(1):36--52, 2008.

\bibitem{thunder}
J.~L. Thunder.
\newblock The number of solutions of bounded height to a system of linear
  equations.
\newblock {\em J. Number Theory}, 43(2):228--250, 1993.

\bibitem{widmer}
M.~Widmer.
\newblock Lipschitz class, narrow class, and counting lattice points.
\newblock {\em Proc. Amer. Math. Soc., to appear}.

\end{thebibliography}
\end{document}